\documentclass[12pt]{article}
\usepackage{amsthm,amsmath,amssymb,latexsym}
\newtheorem{proposition}{Proposition}
\newtheorem{theorem}{Theorem}
\newtheorem*{remark}{Remark}
\newtheorem*{nclemma}{Lemma 5.1 in [Gia]}
\newtheorem*{lemma}{Main Lemma}
\newtheorem*{main}{Theorem $1'$}

\newcommand{\grad}{\nabla}

\begin{document}
\title{Local pointwise  estimates for solutions of the $\sigma_2$ curvature equation on $4$ manifolds}
\author{Zheng-Chao Han \thanks{The author acknowledges the support on this work by  NSF through grant DMS-0103888,
and the encouragement   of Professors H. Brezis, S.-Y. A. Chang and P. Yang.}
\\Department of Mathematics \\Rutgers University \\110 Frelinghuysen Road \\
Piscataway, NJ 08854 \\ \texttt{zchan@math.rutgers.edu}}
\date{}
\maketitle
\begin{abstract}
The study of the $k$-th elementary symmetric function of the Weyl-Schouten curvature tensor of a Riemannian metric,
the so called $\sigma_k$ curvature,
has produced many fruitful results in conformal geometry in recent years, 
especially when the dimension of the underlying manifold is $3$ or $4$. In these studies in conformal 
geometry, the deforming conformal factor is considered to be a solution of  a fully nonlinear elliptic PDE.
Important  advances have been made in recent years
in the understanding of the analytic behavior of solutions
of the PDE,
including the adaptation  of Bernstein type estimates in integral form, global and local derivative estimates, 
classification of entire solutions and analysis of blowing up solutoins. Most of these results require 
derivative bounds on the $\sigma_k$ curvature. The derivative estimates also require an a priori $L^{\infty}$
bound on the solution. This work provides local $L^{\infty}$ and Harnack estimates for
solutions of the $\sigma_2$ curvature equation on $4$ manifolds, under only $L^p$ bounds on the
$\sigma_2$ curvature, and the natural assumption of small volume(or total $\sigma_2$ curvature). 
\end{abstract}

\section{Introduction and Statements of the results}

This paper addresses local $L^{\infty}$ and Harnack estimates for \emph{admissible} solutions $w$ to either
\begin{equation} \label{1}
\sigma_2(g^{-1} \circ A_g) = K(x),
\end{equation}
or
\begin{equation} \label{2}
\sigma_2(g_0^{-1} \circ A_g)  = f(x),
\end{equation}
where $g_0$ is a fixed  background metric on a 4-manifold $M^4$ and $g = e^{2w(x)}g_0$ is 
a metric conformal to $g_0$, $ A_g $ is the the Weyl-Schouten tensor of the  metric $g$, 
\begin{equation} \label{trans}
\begin{split}
A_g &=  \frac {1}{n-2} \{ Ric - \frac {R}{2(n-1)}g \} \\
& = A_{g_0} - \left[ \grad^2w-dw \otimes dw + \frac 12 ||\grad w||^2 g_0\right],
\end{split}
\end{equation}
and $\sigma_k(\Lambda)$, for any $1-1$ tensor $\Lambda$ on an $n-$dimensional  vector space
and $k \in \mathbb N$, $0 \leq k \leq n$, is 
the $k$-th elementary symmetric function of the eigenvalues of $\Lambda$;
$K(x)$ and $f(x)$ are two nonnegative functions with only some \emph{integrability} assumptions, and
an \emph{admissible} solution is defined to be a $w \in C^2(M^4)$ such that for all $x \in M^4$,
$A_g (x) \in \Gamma^+_2$ (see next paragraph for the definition of $\Gamma^+_2$)
 and \eqref{1} or \eqref{2} is satisfied. Note that, since
\[
\sigma_2(g^{-1} \circ A_g) = e^{-4w} \sigma_2(g_0^{-1} \circ A_g),
\]
 so a solution of \eqref{1} is a solution of \eqref{2} with 
$f(x)  = K(x) e^{4w}$.

It is natural to restrict to metrics whose Weyl-Schouten tensor is in  the $\Gamma^{+}_k$ class, 
\emph{i.e.}, those metrics such that $\sigma_j(g^{-1} \circ A_g ) >0$ for $1 \leq j \leq k$, 
 because, \emph{for a metric in such a class with $k >1$}, (i) $\sigma_k(g^{-1} \circ A_g )$ places
a much stronger control on the curvature tensor: Chang, Gursky and Yang 
\cite{CGY1} observed that if 
$\sigma_1(g^{-1} \circ A_g), \sigma_2(g^{-1} \circ A_g) >0$ at a point on a $4$-dimensional manifold, then
the Ricci tensor of $g$ is positive definite at that point; this algebraic relation has been generalized to
higher dimensions by Guan, Viaclovsky, and Wang \cite{GVW}; (ii) the expression  
$\sigma_k (g^{-1} \circ A_{g} )$ is 
 a  \emph{fully nonlinear} PDO in $w$ that becomes \emph{elliptic}. Another reason the study of $\sigma_k$
curvature has attracted strong interest in recent years is its appearance in conforml invariant. Viaclovsky
establishes in \cite{V1} that when $2k =n$ on  a closed manifold $M^n$, and $M$ is locally conformally flat 
if $k>2$, the integral $\int_M \sigma_k (g^{-1} \circ A_g) d\, vol_g$ is conformally invariant. In fact,
in dimension $4$, the $\sigma_2$ curvature comes into play in the 
Chern-Gauss-Bonnet formula. 
 The first important application of the
$\sigma_k$-curvature to conformal geometry is the main theorem in \cite{CGY1}, 
where the authors prove that
if (a) $\int_{M^4} \sigma_2 (A_g) d\,vol_g$, which is conformally invariant on $M^4$, is  positive; and (b)
the Yamabe class of $(M^4, g)$ is positive, then there is a conformal metric $\tilde g = e^{2w}g$ on 
$M^4$ such that $A_{\tilde g} \in \Gamma_2^+$. Note that
 $\sigma_1(A_g)$ is simply a constant multiple of the scalar curvature of $g$, so
 $A_g$ in the $\Gamma^{+}_k$ class is a generalization of the notion that  the scalar
curvature $R_g$ of $g$ having a fixed $+$ sign. Further applications of the $\sigma_k$ curvature
in geometry appear in \cite{CGY2, CGY3, GV1, GV2, GV3}.

Note that $\sigma_k(g^{-1} \circ A_g)$ has a divergence structure as given by
\begin{proposition} \label{divstr}
If $g = e^{2w}|dx|^2$ is locally conformally flat, then
\begin{equation} \label{div}
 k \sigma_k(g^{-1}\circ A_g) = (n-2k) \sum_{j=1}^{k} \frac{\sigma_{k-j}(g^{-1}\circ A_g)}{2^j}||\grad w||_g^{2j}
- \grad_a X^a,
\end{equation}
where
\[
X^a =\left[\sum_{j=1}^{k}\frac{T_{k-j}(g^{-1}\circ A_g)^a_b ||\grad w||_g^{2(j-1)}}{2^{j-1}}\right]\grad^b w,
\]
$T_{k-j}(g^{-1}\circ A_g)^a_b$ is the $(k-j)$-th Newton transform of $g^{-1}\circ A_g$:
\[
T_{k-j}(g^{-1}\circ A_g) = \sum_{i=0}^{k-j}(-1)^i\sigma_{k-j-i}(g^{-1}\circ A_g)(g^{-1}\circ A_g)^i.
\]
In this theorem all norms and differentiation are instrinsic with $g$.
\end{proposition}

We remind the reader that $T_{k-j}(g^{-1}\circ A_g)$ is positive definite for $1 \le j \le k$
 when $A_g \in \Gamma_k^+$, see \cite{V2}.
Also note that when $2k=n$, $\sigma_k(g^{-1}\circ A_g)$ is a pure divergence, confirming the 
results of  \cite{V1} as mentioned above.
In her thesis \cite{G} Gonzalez also exploits the divergence strucuture of $\sigma_k(g^{-1} \circ A_g)$.
See also \cite{G2, G3}.
Her work deals with the case $2k <n$ as opposed to our case where $2k=n$. The case she deals with
 exhibits somewhat different analytical behavior, as hinted
by the extra terms in \eqref{div} with definite signs after $(n-2k)$.
She does not give an expression such as \eqref{div}, instead, uses an inductive relation between
$\sigma_k(g^{-1} \circ A_g)$ and $\sigma_{k-1}(g^{-1} \circ A_g)$. \eqref{div} appears to be based on 
the same analytical structure as exploited in \cite{G}.

%The $k$-th Newton transform of $\Lambda$ is defined as 
%\[
%T_k(\Lambda) = \sum_{j=0}^{n}(-1)^j\sigma_{k-j}(\Lambda)\Lambda^j.
%\]
%The following elementary algebraic properties will be useful later on.
%

%\begin{proposition} 
%\begin{enumerate}
%\item[(i)] $(k+1)\sigma_k(\Lambda) = T_k (\Lambda)^a_b  \Lambda^b_a$.
%\item[(ii)]  $(n-k)\sigma_k(\Lambda) = T_k (\Lambda)^a_a$.
%\item[(iii)] $T_k (\Lambda)^a_b = \sigma_k(\Lambda) \delta^a_b - T_{k-1}(\Lambda)^a_c  \Lambda^c_b$.
%\item[(iv)] If $\Lambda \in \Gamma^+_k$, \emph{i.e.}, $\sigma_j(\Lambda) >0$, for $1\le j \le k$, then
%$T_{k-j}(\Lambda)$ is positive definite, for $1\le j \le k$.
%\end{enumerate}
%\end{proposition}
%$\Lambda = g^{-1} \circ A_g$ as a $1-1$ \emph{tensor field} enjoys  the following property.
%
%\begin{proposition}
%If $g$ is locally conformally flat, or if $k = 1$, then $T_k (g^{-1} \circ A_g)^a_{b,a} = 0$.
%\end{proposition}

 Since most of the geometric applications involving $\sigma_k$ deal with the case of $k=2$ in
dimension 4, we limit ourselves to this case in this paper.
For simplicity of presentation, we will take the background metric
$g_0 = |dx|^2$ to be the flat metric in a ball in $\mathbb R^4$.
 For a general $g_0$, only minor computational changes are needed.
In our case, the left hand sides of   \eqref{2} can be written as divergence in the background
metric $g_0$:
\begin{equation} \label{ibp}
2\sigma_2 (g_0^{-1} \circ A_g) = - \partial_a(M^a_b \partial^b w),
\end{equation}
where
\[
M^a_b =  T_1(g_0^{-1} \circ A_g)^a_b + \frac {|\grad w|^2}{2} \delta^a_b.
\]
Here and in the remaining of the paper, $\grad w$ denotes the gradient of $w$ in the background metric $g_0$.
\eqref{ibp}, which follows from \eqref{div}, is already exploited in \cite{CGY2}.

It turns out to make sense to discuss a weak notion of admissible solution (subsolution, supersolution) of
\eqref{1} or \eqref{2}. An admissible  $W^{2,2}$ 
solution (subsolution, supersoluton) of \eqref{1} or \eqref{2} is a $w \in W^{2,2}(M^4)$ such that, 
for \emph{a.e.} $x$, $A_g (x) \in \Gamma^+_2$, and the left hand side of 
\eqref{1} or \eqref{2} $=(\leq ,\,\geq)\;$ the right hand side.
In Propositions~\ref{p>1} and \ref{whnk} we will consider pointwise upper (lower) bound and weak Harnack 
inequality for a $W^{2,2}$ admissible subsolution (supersolution) of \eqref{2}, respectively. In 
Theorems~\ref{hnkforw} and \ref{hnk} we provide Harnack inequality for $W^{2,2}$ 
admissible solutions of
\eqref{2}. For geometric applications involving solutions of \eqref{1},
the following estimate under small volume is probably most useful.
\begin{theorem} \label{small} Assume $0 \leq \inf_{B_{2R}} K \leq \sup_{B_{2R}} K < \infty$. 
There exist absolute  constant $\epsilon_0 >0$ small and $C^*>0$ depending on
$(\sup_{B_{2R}} K )\int_{B_{2R}} e^{4w} \, d vol_{g_0}$
such that for any admissible solution $w$ of \eqref{1}  on 
$B_{2R} \subset \mathbb R^4$, 
if 
\begin{equation}\label{se}
 \int_{B_{2R}}K(x) e^{4w} \, d vol_{g_0} < \epsilon_0,
\end{equation} 
then
\[
\sup_{B_R} e^{w} \leq C^* \left(\frac{1}{R^4} \int_{B_{2R}} e^{4w} \, d vol_{g_0}\right)^{\frac 14},
\]
and
\[
\sup_{B_R} e^{w} \leq C^* \inf_{B_R} e^{w}.
\]
\end{theorem}
\begin{remark}
From the proof it will be clear that the condition $\sup_{B_{2R}} K < \infty$ may be relaxed to
$||K||_p < \infty$ for some $p>1$. Then $\epsilon_0$  depends on $p$ and $C^*$ also depends on $p$
as well as on $R^{-4/p} ||K||_{L^p(B_{2R})} \int_{B_{2R}} e^{4w} \, d vol_{g_0}$.
If one is willing to assume $\sup_{B_{2R}} K < \infty$ and  the smallness of 
$(\sup_{B_{2R}} K) \int_{B_{2R}} e^{4w} \, d vol_{g_0}$, then 
the following slightly different, less geometric version of Theorem~\ref{small}, is much easier to prove.
\end{remark}
\begin{main} 
Assume  $0 \leq \inf_{B_{2R}} K \leq \sup_{B_{2R}} K < \infty$. 
There exist absolute  constants $\epsilon_1 >0$ small and $C_*>0$
such that for any admissible solution $w$ of \eqref{1}  on 
$B_{4R} \subset \mathbb R^4$, 
if 
$$
(\sup_{B_{2R}} K )\int_{B_{2R}} e^{4w} \, d vol_{g_0} < \epsilon_1,
$$ 
then
\[
\sup_{B_R} e^{w} \leq C_* \left(\frac{1}{R^4} \int_{B_{2R}} e^{4w} \, d vol_{g_0}\right)^{\frac 14},
\]
and
\[
\sup_{B_R} e^{w} \leq C^* \inf_{B_R} e^{w}.
\]
\end{main}

Theorems~\ref{small} and $1'$
are consequences of  the following propositions and theorems, which provide  pointwise and 
Harnack estimates for
solutions/subsolutions/supersolutions of the fully nonlinear equation \eqref{2}, and are of independent interest. 
\begin{proposition} \label{p>1}
Assume $w$ is an admissible $W^{2,2}$ subsolution of \eqref{2} on $B_{2R} \subset \mathbb R^4$,
 and $||f||_{L^p(B_{2R})} < \infty$ for some $p>1$. Then $w$ is bounded from above on $B_R$ and for any $\beta >0$,
 there exists 
\[
C=C(p,R^{4(1-\frac {1}{p})}||f||_{L^p(B_{2R})}, \beta)>0
\]
 such that 
\begin{equation} \label{upper}
\sup_{B_R} e^{w} \leq C \left(\frac{1}{R^4} \int_{B_{2R}} e^{\beta w} d\, vol_{g_0} \right)^{\frac {1}{\beta}},
\end{equation}
and there exists $C^*=C^*(p, \beta)>0$ such that,
if $m \leq w \leq M$ on $B_{2R}$, and
$\gamma = \max (1, R^{\frac {4}{3} (1 -\frac {1}{p})} ||f||^{\frac {1}{3}}_{L^p(B_{2R})})$, then
\begin{equation} \label{supwm}
\sup_{B_R} (\gamma +  w-m) \leq C^* \left( \frac{1}{R^4} \int_{B_{2R}} (\gamma + w-m)^{\beta} \right)^{\frac {1}{\beta}},
\end{equation}
\begin{equation} \label{whnkMw}
\left(\frac{1}{R^4} \int_{B_{2R}} (\gamma + M - w)^{\beta} \right)^{\frac {1}{\beta}} \leq C^* \inf_{B_R} (\gamma + M - w).
\end{equation}
\end{proposition}
Here and in the following, the integrals are taken with respect to the measure generated by $g_0$.
\begin{proposition} \label{whnk}
Assume $w$ is an admissible $W^{2,2}$ supersolution of \eqref{2} on $B_{2R} \subset \mathbb R^4$,
 and $||f||_{L^p(B_{2R})} < \infty$ for some $p>1$. Then $w$ is bounded from below on $B_R$ and for any $\beta >0$,
 there exists 
\[
C=C(p,R^{4(1-\frac {1}{p})}||f||_{L^p(B_{2R})}, \beta)>0
\]
such that
\begin{equation} \label{lower}
\inf_{B_R} e^{w} \geq C \left(\frac{1}{R^4} \int_{B_{2R}} e^{-\beta w} d\, vol_{g_0} \right)^{\frac {1}{-\beta}},
\end{equation}
and there exists $C^*=C^*(p, \beta)>0$ such that,
if $m \leq w \leq M$ on $B_{2R}$, and
$\gamma = \max (1, R^{\frac {4}{3} (1 -\frac {1}{p})} ||f||^{\frac {1}{3}}_{L^p(B_{2R})})$, then
\begin{equation} \label{supMw}
\sup_{B_R} (\gamma + M - w) \leq C^* \left( \frac{1}{R^4} \int_{B_{2R}} (\gamma + M - w)^{\beta} \right)^{\frac {1}{\beta}},
\end{equation}
\begin{equation} \label{whnkwm}
\left( \frac{1}{R^4} \int_{B_{2R}} (\gamma + w-m)^{\beta} \right)^{\frac {1}{\beta}} \leq C^* \inf_{B_R} (\gamma + w-m).
\end{equation}
\end{proposition}
\begin{remark}
As will be seen later in the proofs, the estimates of Proposition~\ref{whnk} follow from the same scheme of
proof as for Proposition~\ref{p>1}. However, slightly different formulations of the estimates of 
Proposition~\ref{whnk} follow
easily by the superharmonicity of $w$, which is a consequence of $A_{e^{2w}g_0} \in \Gamma_2^+$. We formulate them in conjuction
with those of Proposition~\ref{p>1}, because, together, they give the following
\end{remark}
\begin{theorem} \label{hnkforw}
Let $w$ be an admissible $W^{2,2}$ solution of \eqref{2} on $B_{2R} \subset \mathbb R^4$,
 and $||f||_{L^p(B_{2R})} < \infty$ for some $p>1$. 
Set $\gamma = \max (1, R^{\frac {4}{3} (1 -\frac {1}{p})} ||f||^{\frac {1}{3}}_{L^p(B_{2R})})$.
Then there exists $C^*=C^*(p)>0$ such that,
\begin{enumerate}
\item[(i)] if  $M$ is   an upper bound of $w$ on $B_{2R}$, then we have
\begin{equation} \label{hnkMw}
\sup_{B_R} (\gamma + M -w) \leq C^* \inf_{B_R} (\gamma + M -w),
\end{equation}
\item[(ii)] if $m$ is a lower bound of $w$ on $B_{2R}$, then we have
\begin{equation} \label{hnkwm}
\sup_{B_R} (\gamma + w-m) \leq C^* \inf_{B_R} (\gamma + w-m).
\end{equation}
\end{enumerate}
\end{theorem}
%\begin{theorem} \label{hld}
%Any admissible $W^{2,2}$ solution $w$ of \eqref{1} with $K \in L^{p}$ for some $p>1$, or of \eqref{2} with $f \in L%^p$ 
%for some $p>1$ is locally  H\"older continuous.
%\end{theorem}
A different formulation of Harnack estimate in terms of $e^w$ is given as
\begin{theorem} \label{hnk}
If $w$ is an admissible solution of \eqref{2} on $B_{2R} \subset \mathbb R^4$,
 and $||f||_{L^p(B_{2R})} < \infty$ for some $p>1$, then there exists 
$C=C(p,R^{4(1 -\frac {1}{p})}||f||_{L^p(B_{2R})})$ such
that
\begin{equation} \label{hnkew}
\sup_{B_R} e^{w} \leq C \inf_{B_R} e^{w}.
\end{equation}
\end{theorem}
Let us put our results into perspective. Of the many important, recent analytic contributions on related 
problems, this work is more directly related to \cite{V2, GW1}, although it is also closely
related to \cite{CGY1}-\cite{CGY2}
and \cite{G}.  For equations of the same type as \eqref{1}, but
with general $k$ and $n$,   Viaclovsky establishes in \cite{V2} \emph{global} $C^1$ and $C^2$ estimates
for $C^4$ admissible solutions, \emph{assuming  $C^0$ estimate on the solution}.
In \cite{CGY1}-\cite{CGY2}, Chang, Gursky and Yang develop important integral estimates for
related equations, some for
a singularly perturbed fourth order equation. 
Later  Guan and Wang
establish in \cite{GW1} \emph{local} $C^1$ and $C^2$ estimates for 
$C^4$ admissible solutions of  equations similar to those in \cite{V2},  
\emph{assuming a one-sided  $C^0$ estimate on the solution}. 
All these results
require \emph{derivative bounds of $K(x)$}, of course, also provide stronger estimates, namely
derivative estimates. Similar results were then proved for a more
 general class of fully nonlinear equations by A. Li and YanYan Li in \cite{LL1}-\cite{LL6}, 
where they also establish Liouville type theorems,
Harnack type theorems in the sense of Schoen, compactness and existence results. 
All these are very important and useful results. In fact, if $C^1$ bounds on $K$ are allowed, then
the conclusion of Theorem~\ref{small} is covered by \cite{GW1}.
However, for some applications in blowing up analysis, the
derivative bounds of $K(x)$ in the assumptions of \cite{V2} and \cite{GW1}
 are  absent, and only some $L^p$ bounds are under control. Our Theorem~\ref{small} 
provides a partial substitute. 

Using our results, we establish the following theorem in a joint work with S.-Y. A. Chang and P. Yang \cite{CHY1}.
\begin{theorem} \label{chy}
Consider a family of admissible conformal metrics $g_j=e^{2w_j}g_c$ on $\mathbb S^4$ with
$\sigma_2(g^{-1}_j\circ A_{g_j}) = K(x)$, where $g_c$ denotes
the canonical round metric on $\mathbb S^4$ and $K(x)$ denotes a fixed $L^{\infty}$ function on $\mathbb S^4$
with a positive lower bound.
\begin{enumerate}
\item There exists at most one isolated simple blow up point in the sense that, 
if $\max w_j = w_j (P_j) \to \infty$, then there exists conformal automorphism $\varphi_j$ of $\mathbb S^4$
such that, if we define $v_j (P) = w_j \circ \varphi_j (P) + \ln | \det (d \varphi_j)|$, we have
\[
v_j (P) - \frac {1}{4} \ln \frac {6}{K(P_j)} \to 0 \quad \text{in} \quad L^{\infty} \quad \text{and} \quad
\int_{\mathbb S^4} |\grad v_j|^4 \to 0.
\]
\item If, furthermore, $K(x)$ is $C^2$ and satisfies a non-degeneracy condition
\[
\Delta K(P) \neq 0 \quad \text{whenever}\quad \grad K(P) =0.
\]
Then there exists apriori $C^{2, \alpha}$ estimate on $w_j$ depending on $\max K$, $\min K$, the $C^2$ norm of
$K$ and the modulus of continuity of $\grad^2 K$.
\end{enumerate}
\end{theorem}
Propositions~\ref{p>1}, \ref{whnk}, Theorems~\ref{hnkforw} and \ref{hnk}
are proved by a Moser iteration scheme. However, the fully nonlinear equations \eqref{1} or \eqref{2},
when regarded as an elliptic equation in $w$ in divergence form, is \emph{not} uniformly elliptic. Moser iteration
procedures have been successfully employed to deal with some non-uniformly elliptic quasilinear equations, see
the classic paper \cite{S}, and, for general reference,  also the monographs 
 \cite{GT, LU}. But they do not directly apply to our situation.
We establish the Moser iteration scheme by exploiting the special divergence structure in the equation through the 
following Lemma.  
\begin{lemma} \label{ml}
Let $G(w)$ be a nonnegative Lipschitz function of $w$. If $w$ is an admissible subsolution of \eqref{2}, 
we will require
$G'(w) \geq 0$; and if $w$ is an admissible supersolution of \eqref{2}, we will require $G'(w) \leq 0$.
Let $\eta \in C^2_0(B_2)$ be a non-negative
cut-off function
on $B_2$ satisfying $\eta \equiv 1$ on $B_1$ and $\eta |\grad^2 \eta| \lesssim |\grad \eta|^2$.  Then
\begin{equation} \label{test}
\int_{B_2} \eta^4 |G'(w)| |\grad w|^4  \lesssim \int_{B_2} \eta^2 |G(w)| \left[|\grad w|^2 |\grad \eta|^2 +
|\grad w|^3 \eta |\grad \eta|  \right] + \int_{B_2} \eta^4 |G(w)| |f|,
\end{equation}
here and in the following, the integrations are all done with respect to the background metric $g_0$,
and we write $X \lesssim Y$ when there is an absolute constant $c>0$, depending perhaps only on dimension,
such that $X \leq c Y$.
\end{lemma}

It turns out that \cite{CGY2} already used an integral estimate like the one in the Main Lemma for a specific 
$f \equiv 0$ and $G(w)=w-\bar w$, where $\bar w$ is the average of $w$ on $B_{2R}$, although they only used that
as a step in the classification of entire solutions and did not pursue 
the iteration of such integral estimates as done here.
Most of the results here were obtained in 2002; some were in slightly different formulations and had different
proofs. In particular, my orginal formulation and
proof of Theorem 1, using harmonic approximation, requied some boundary
information of the solution. I wish to thank Professors Chang and Yang for their questioning of my earlier
proof, which prompted me to find the current better proof.
%The initial proofs were along similar lines as described here; the current proofs are more clean. 
I would also like to call attention to  \cite{G, G2, G3}, 
where, as mentioned earlier,
M. Gonzalez  also exploits the divergence structure of $\sigma_k$ and adapts the Moser iteration scheme.
However, other than these two similarities on methods, which were developed independently --- we didn't 
learn of each other's work until after we both completed our work and began to report on them, there is 
no overlap between our work.
 Our work originated from different motivations and address different situations: this work
was motived mainly for applications in apriori estimates such as in Theorem~\ref{chy}; while Gonzalez's
work was mostly for studying the  size of the singular set of solutions to the $\sigma_k$ 
equations of the type similar to \eqref{1} with $2k < n$ and their removability. In fact, one of Gonzalez's theorems
says that, when $2k<n$, 
 an isolated singularity of the $\sigma_k$ Yamabe equation with finite volume is removable. The corresponding statement
in our case does not hold, despite our local $L^{\infty}$ estimates under small volume. This 
 can be seen from solutions in my joint work \cite{CHY2} with S.-Y. A. Chang and P. Yang. 

\section{Sketch of proofs}
We will first indicate  how Theorem $1'$ follows from Proposition~\ref{p>1} and Theorem~\ref{hnk}.
Then we will provide a proof for the Main Lemma, which is the basis for all the iteration procedures. Finally
we will describe the proofs for  Propositions~\ref{p>1}, \ref{whnk}, Theorems~\ref{hnkforw}, \ref{hnk}, and \ref{small}.
The proof for Proposition 1 is in fact quite routine, making use of the transformation
formulas such as (3) and the fact that $\grad_a T_i (g^{-1} \circ A_g)^a_b =0$ in our situation. Since
it is not used essentially in this work, it will be omitted here and will be provided in a future work.

\begin{proof}[Proof of Theorem $1'$] Here we can take $R$ to be $1/2$.
 The general case follows from rescaling.
\eqref{1} has translation covaraince: if we set $\widetilde w = w + \frac {1}{4} \ln ||K||_{\infty}$, then
$\widetilde w $ satisfies
\[
\sigma_2 (g_0^{-1} \circ A_{\widetilde w } ) = \widetilde K (x) e^{4\widetilde w },
\]
where $ \widetilde K (x)= K(x) /||K||_{\infty}$.
It is easier to work with $u(x) = e^{\widetilde w(x)}$. For the local upper bound, 
we only need to bound $\sup_{|x| \leq 1} (1-|x|) u(x)$ in terms
of $\int_{B_1} u^4(x) dx$. Note that $\sup_{|x| \leq 1} (1-|x|) u(x) < \infty$ by Proposition~\ref{p>1}.
Let $x_0$ with $|x_0| <1$ satisfy
\[
\sup_{|x| \leq 1} (1-|x|) u(x) \leq 2 (1-|x_0|)u(x_0).
\]
Set $1-|x_0| = 2r_0$
and $v(z) = \rho u(x_0 + \rho z)$, with $\rho >0$ chosen so that $v(0)=1$. When $|z| \leq r_0/\rho$, we have
$1-|x_0+\rho z| \geq r_0$, so that
\[
r_0u(x_0+\rho z) \leq (1-|x_0+\rho z|) u(x_0+\rho z) \leq 4r_0 u(x_0).
\]
Thus $v(z) \leq 4$ for $|z| \leq r_0/\rho$. Note that $u^2(x)|dx|^2 = v^2(z)|dz|^2$, so that
\[
\sigma_2(g_0^{-1} \circ A_v)=\widetilde K(x_0 + \rho z) v^4,
\]
where $A_v$ is the Weyl-Schouten tensor of $v^2(z)|dz|^2$. If $r_0/\rho \geq 1$, then $v(z)\leq 4$ on $|z| \leq 1$.
The conditions for Proposition~\ref{p>1}
are satisfied on $|z| \leq 1$. Noting that $|| \widetilde K||_{\infty} =1$.
So we apply Proposition~\ref{p>1} with $R = 1/2$ and $p=\beta =4$ to obtain
an absolute constant $C_* >0$ such that
\[
1 =v(0) \leq C_* \left( \int_{|z| \leq 1} v^4 (z) dz\right)^{\frac {1}{4}}= C_* 
\left( \int_{B(x_0, r_0)} u^4(x) dx \right)^{\frac {1}{4}}\leq C_*
\left( || K||_{\infty}  \int_{B_1} e^{4w(x)} dx \right)^{\frac {1}{4}}.
\]
This can't happen if $|| K||_{\infty} \int_{B_1} e^{4w(x)} dx \leq \epsilon_1$ and $\epsilon_1 C_*^4 < 1$.
 Choose and fix such
an $\epsilon_1$. Then we must have $r_0/\rho < 1$. Again we can apply Proposition~\ref{p>1} on
$|z| \leq r_0/\rho$ to obtain
\[
1 \leq C_*^4 \left( \frac {\rho}{r_0}\right)^4 \int_{|z| \leq \frac {r_0}{\rho}} v^4 (z) dz
\leq C_*^4 \left( \frac {\rho}{r_0}\right)^4 \int_{|x| \leq 1} u^4(x)dx.
\]
Recall that $\rho u(x_0) =1$. So we have
\[
r_0^4 u(x_0)^4 \leq C_*^4 \int_{|x| \leq 1} u^4(x)dx,
\]
from which it follows that
\[
\sup_{|x| \leq 1} (1-|x|) u(x) \leq 4r_0 u(x_0) \leq 4 C_* \left( \int_{|x| \leq 1} u^4(x)dx \right)^{\frac {1}{4}}.
\]
For the Harnack estimate, we will apply Theorem~\ref{hnk} with $f = K(x) e^{4w}$. For that purpose we will need
to have an upper bound on $R^{4(1-\frac 1p)}||f||_{L^p(B_{2R})}$. An almost identical verification is carried out
in the proof of Theorem 1 later. Please refer to that part of the proof.
\end{proof}

\begin{proof}[Proof of the Main Lemma] Recall that for $g= e^{2w}g_0 = e^{2w} |dx|^2$,
\[
 2\sigma_2 (g_0^{-1} \circ A_g) = - \partial_a(M^a_b \partial^b w),
\]
where
\[
M^a_b =  T_1(g_0^{-1} \circ A_g)^a_b + \frac {|\grad w|^2}{2} \delta^a_b.
\]
We obtain
\[
\begin{split}
\int_{B_2} 2 \sigma_2 (g_0^{-1} \circ A_g) \eta^4 G(w) & = \int_{B_2}  M^a_b \partial_a w \partial^b\left[ \eta^4 G(w) \right] \\
		&= \int_{B_2} \eta^4 G'(w) M^a_b \partial_a w \partial^b w  + 
				\int_{B_2} 4 \eta^3 G(w) M^a_b \partial_a w \partial^b \eta \\
	&\begin{cases}
		\geq \frac 12 \int_{B_2} |\grad w|^4 \eta^4 G'(w) + 
			\int_{B_2} 4 \eta^3 G(w) M^a_b \partial_a w \partial^b \eta , &\text{if $G' \geq 0$;} \\
                \leq \frac 12 \int_{B_2} |\grad w|^4 \eta^4 G'(w) + 
			\int_{B_2} 4 \eta^3 G(w) M^a_b \partial_a w \partial^b \eta , &\text{if $G' \leq 0$.}
	\end{cases}
\end{split}
\]
Here and in the following of the proof, all the integration by parts used  can be justified for $W^{2,2}$ 
admissible solutions. Thus, in all cases, we have
\[
\int_{B_2} |\grad w|^4 \eta^4 |G'(w)| \leq 8| \int_{B_2}  \eta^3 G(w) M^a_b \partial_a w \partial^b \eta| + 
		4 \int_{B_2} |f(x)| \eta^4 |G(w)|.
\]
$M^a_b$  is an expression involving the first and
second derivatives of $w$, so, apriori, we have \emph{no upper bound} 
on the eigenvalues of $M^a_b$. However, it follows from \eqref{trans} and the definition of $M^a_b$ that
\[
M^a_b = w^a_b - (\Delta w) \delta^a_b -w^a w_b,
\]
and
\[
\begin{split}
2 M^a_b \partial_a w &= \partial_b(|\grad w|^2) - 2(\Delta w) \partial_b w - 2 |\grad w|^2 \partial_b w \\
		& = 2 \partial_b(|\grad w|^2) - \partial_a (2 \partial^a w \partial_b w) 
				-2 |\grad w|^2 \partial_b w ,
\end{split}
\]
We can integrate by parts to estimate the integral 
\[
\begin{split}
\int_{B_2}  \eta^3 G(w) M^a_b \partial_a w \partial^b \eta 
= & \int_{B_2} \left(\partial^a w \partial_b w - |\grad w|^2 \delta^a_b \right)
             \partial_a \left( \eta^3 \partial^b \eta G(w) \right) - 
             \eta^3 G(w) |\grad w|^2 \partial_b w \partial^b \eta  \\
= & \int_{B_2} \left(\partial^a w \partial_b w - |\grad w|^2 \delta^a_b \right)
             \{ \left( 3 \eta^2 \partial_a \eta \partial^b \eta + \eta^3 \partial_a^b \eta \right) G(w) + \\
	& +G'(w) \eta^3 \partial_a w \partial^b \eta \} 
    -   \int_{B_2}   \eta^3 G(w) |\grad w|^2 \partial_b w \partial^b \eta \\
= &  \int_{B_2} \left(\partial^a w \partial_b w - |\grad w|^2 \delta^a_b \right)
              \left( 3 \eta^2 \partial_a \eta \partial^b \eta + \eta^3 \partial_a^b \eta \right) G(w) \\ 
	&- \int_{B_2}	\eta^3 G(w) |\grad w|^2 \partial_b w \partial^b \eta .
\end{split}
\]
Using $|\eta \partial_a^b \eta| \lesssim |\grad \eta|^2$, we conclude the proof of the Main Lemma.
\end{proof}

\begin{proof}[Proof of Propositions~\ref{p>1} and  \ref{whnk}] 
\cite{S} serves as a useful guide in the adaptation of Moser's iteration procedure to our situation.
The proof of \eqref{upper} and \eqref{lower}
is done by plugging in \eqref{test} $G(w) = e^{4\beta w}$, or
more strcitly speaking, truncations of $e^{4\beta w}$ at large $|w|$.  For simplicity, we will not do the truncation
in the test function; Instead, we will demonstrate the  estimates apriori by
simply plugging  $G(w) = e^{4\beta w}$ in \eqref{test}. Then 
\begin{equation}
	\begin{split}
	\int_{B_2} \eta^4 |G'(w)| |\grad w|^4 &= 4|\beta| \int_{B_2} \eta^4 e^{4\beta w} |\grad  w|^4\\
					& = 4|\beta|^{-3} \int_{B_2} \eta^4 |\grad e^{\beta w}|^4,
	\end{split} 
\end{equation}
\begin{equation}
	\begin{split}
	\int_{B_2} \eta^2 |G(w)| |\grad w|^2 |\grad \eta|^2 &= \beta^{-2} \int_{B_2} \eta^2 |\grad e^{\beta w }|^2
			 e^{2\beta w} |\grad \eta|^2 \\
			 & \leq \beta^{-2} \left(\int_{B_2} \eta^4 |\grad e^{\beta w}|^4\right)^{\frac 12}
		\left(\int_{B_2} e^{4\beta w} |\grad \eta|^4 \right)^{\frac 12}\\
			&\leq \frac {1}{2|\beta|^{3}} \int_{B_2} \eta^4 |\grad e^{\beta w}|^4 + 
				\frac {1}{2|\beta|} \int_{B_2} e^{4\beta w} |\grad \eta|^4 ,
	\end{split}
\end{equation}
and
\begin{equation}
	\begin{split}
	\int_{B_2} \eta^3  |G(w)||\grad w|^3|\grad \eta| &= |\beta|^{-3} \int_{B_2} \eta^3 |\grad e^{\beta w}|^3 
			e^{\beta w }|\grad \eta| \\
		&\leq |\beta|^{-3} \left(\int_{B_2} \eta^4 |\grad e^{\beta w}|^4\right)^{\frac 34}
				\left(\int_{B_2} e^{4\beta w} |\grad \eta|^4 \right)^{\frac 14}\\
                                &\leq \frac {3}{4|\beta|^3} \int_{B_2} \eta^4 |\grad e^{\beta w}|^4 + 
					\frac {1}{4|\beta|^3}  \int_{B_2} e^{4\beta w} |\grad \eta|^4.
	\end{split}
\end{equation}
Combining the above, we have
\begin{equation} \label{pfinal}
 \int_{B_2} |\grad (\eta e^{\beta w})|^4 \lesssim (1+\beta^{2}) \int_{B_2} e^{4\beta w} (|\grad \eta|^4 +  \eta^4)
		+ |\beta|^3 ||f||_p ||\eta e^{\beta w}||_{4p'}^4,
\end{equation}
where $p' = \frac {p}{p-1}$ is the H\"older conjugate of $p$. Set $q = 8 p'$ and $\theta = (2p'-1)^{-1}$. Then
$0 < \theta <1$ and 
\[
\frac {1}{4p'} = \frac {1-\theta}{q} + \frac {\theta}{4}.
\]
We use interpolation to estimate the last term in \eqref{pfinal},
$||\eta e^{\beta w}||_{4p'}$, in terms of $||\eta e^{\beta w}||_4$ and
$||\eta e^{\beta w}||_q$:
\begin{equation} \label{hder}
||\eta e^{\beta w}||_{4p'}^4 \leq ||\eta e^{\beta w}||_4^{4\theta} ||\eta e^{\beta w}||_q^{4(1-\theta)},
\end{equation}
and use the  Sobolev inequality to estimate $||\eta e^{\beta w}||_q$.  Putting these in 
\eqref{pfinal}, we have
\[
\begin{split}
||\grad \left(\eta e^{\beta w}\right)||_4^{4} & \lesssim (1+\beta^{2}) ||(|\grad \eta| +  \eta) e^{\beta w}||^4_4 
 + q^{3(1-\theta)}|\beta|^3 ||f||_p ||\eta e^{\beta w}||_4^{4\theta} ||\grad \left(\eta e^{\beta w}\right)||_4^{4(1-\theta)} \\
 & \lesssim (1+\beta^{2}) ||(|\grad \eta| +  \eta) e^{\beta w}||^4_4
+ \theta \left(|\beta| q^{1-\theta} ||f||_p^{\frac 13}\right)^{\frac 3 \theta} ||\eta e^{\beta w}||_4^{4} +
(1-\theta) ||\grad \left(\eta e^{\beta w}\right)||_4^{4}.
%q^3 ||\grad (\eta e^{\beta w})||_4^4 \\
%	&\lesssim q^3 (1+\beta^{2}) ||(|\grad \eta| + \eta) e^{\beta w}||_4^4 + 
%	(|\beta| q)^3 ||f||_p ||\eta e^{\beta w}||_4^{4\theta} ||\eta e^{\beta w}||_q^{4(1-\theta)} \\
%	& \lesssim (1-\theta) ||\eta e^{\beta w}||_q^{4} + 
%			\theta (|\beta| q)^{3/\theta}||f||_p^{1/\theta} ||\eta e^{\beta w}||_4^{4} +
%			q^3 (1+\beta^{2}) ||(|\grad \eta| + \eta) e^{\beta w}||_4^4 .
\end{split}
\]
Thus
\[
||\grad \left(\eta e^{\beta w}\right)||_4^{4} \lesssim  
\left(|\beta| q^{1-\theta} ||f||_p^{\frac 13}\right)^{\frac 3 \theta} ||\eta e^{\beta w}||_4^{4} + \theta^{-1} 
(1+\beta^{2}) ||(|\grad \eta| + \eta) e^{\beta w}||_4^4 ,
\]
and
\[
||\eta e^{\beta w}||_q^{4}  \lesssim (|\beta| q)^{3/\theta}||f||_p^{1/\theta} ||\eta e^{\beta w}||_4^{4} +
					\theta^{-1} q^3 (1+\beta^{2}) ||(|\grad \eta| + \eta) e^{\beta w}||_4^4.
\]
Since $\frac {3}{\theta} = 3(2p'-1) > 2$, we can write the above estimate as
\begin{equation} \label{final}
||\eta e^{\beta w}||_q \leq M (1+|\beta|)^{\frac {3}{4\theta}} ||(|\grad \eta| + \eta) e^{\beta w}||_4,
\end{equation}
where $M$ is a constant depending on $||f||_{p}$ and $p$:
\[
M \sim \left[p' ||f||_{p} \right]^{(2p'-1)/4} + p'.
\]
From \eqref{final}, one can invoke the Moser iteration procedure to prove \eqref{upper} and \eqref{lower}.

To prove \eqref{whnkMw} and \eqref{supMw}, for instance, we first
plug $G(w) = (\gamma + M -w)^{4\beta-3}, \beta \neq 0, \frac 34$  
into \eqref{test}. Then
\[
\eta^4 |G'(w)| |\grad w |^4 = |4\beta-3|\beta^{-4} \eta^4 |\grad (\gamma + M-w)^{\beta}|^4,
\]
\[
\begin{split}
\eta^2 |G(w)| |\grad w|^2 |\grad \eta|^2 &= \beta^{-2} \eta^2 |\grad (\gamma + M-w)^{\beta}|^2 (\gamma + M-w)^{2\beta-1} |\grad \eta|^2 \\
&\leq \frac{\epsilon^2}{2\beta^4} \eta^4 |\grad (\gamma + M-w)^{\beta}|^4 + \frac {1}{2\epsilon^2} (\gamma + M-w)^{4\beta-2}|\grad \eta|^4 \\
&  \leq \frac{\epsilon^2}{2\beta^4} \eta^4 |\grad (\gamma + M-w)^{\beta}|^4 + \frac {\gamma^{-2}}{2\epsilon^2} (\gamma + M-w)^{4\beta} |\grad \eta|^4,
\end{split}
\]
\[
\begin{split}
\eta^3 |\grad w|^3 |G(w)| | \grad \eta| &= |\beta|^{-3} \eta^3  |\grad (\gamma + M-w)^{\beta}|^3 (\gamma + M-w)^{\beta} |\grad \eta| \\
&\leq \frac {3 \epsilon^{4/3}}{4\beta^4} \eta^4 |\grad (\gamma + M-w)^{\beta}|^4 + \frac {1}{4\epsilon^4} (\gamma + M-w)^{4\beta}|\grad \eta|^4, 
\end{split}
\]
and
\[
\eta^4 |G(w)| |f| = \eta^4 (\gamma + M-w)^{4\beta-3} |f| \leq \gamma^{-3} \eta^4 |f| (\gamma + M-w)^{4\beta}.
\]
Choosing $\epsilon >0$ small (its size can be independent of $\beta$ as long as $4\beta -3$ stays away from $0$),
and treating $\int_{B_2} \gamma^{-3} \eta^4 |f| (\gamma + M-w)^{4\beta}$ as in \eqref{pfinal} and \eqref{hder}---noting that $\gamma = \max (1, ||f||_p^{1/3})$, we have
\begin{equation} \label{poweriteration}
|4\beta-3| \int_{B_2}  \eta^4 |\grad (\gamma + M-w)^{\beta}|^4 \lesssim \beta^4  \int_{B_2} (\gamma + M-w)^{4\beta} (|\grad \eta|^4 + \eta^4).
\end{equation}
To complete the proof of \eqref{supMw}, we take $\beta > 3/4$. Then $G'(w) \leq 0$, so we  have
\eqref{poweriteration} and can use it to carry out the Moser iteration to obtain \eqref{supMw} with the restriction
$\beta > 3$ there. But this restriction can be dropped via a device such as Lemma 5.1 in \cite{Gia}.
\begin{nclemma}
Let $E(t)$ be a nonnegative bounded function on $0\leq T_0 \leq  T_1$. Suppose
that for $T_0 \leq t < s \leq T_1$ we have
\[
E(t) \leq \theta E(s)+ A (s-t)^{-\alpha} + B 
\]
with $\alpha >0$, $0\leq \theta < 1$ and $A, B$ nonnegative constants. Then there exists a 
constant $c = c(\alpha, \theta)$, such that for all $T_0 \leq \rho < R \leq T_1$ we have
\[
E(\rho) \leq c \left[ A (R-\rho)^{-\alpha} +B \right]
\]
\end{nclemma}
To complete the proof of \eqref{whnkMw},  we can use \eqref{poweriteration} with $\beta < 3/4, \neq 0$
and need to appeal to the 
John-Nirenberg Theorem. For this purpose, we plug
$G(w) = (\gamma + M -w)^{-3}$ into \eqref{test}. Similar computations as above lead to
\begin{equation} \label{jn}
\int_{B_{2R}} \eta^4 | \grad \log ( \gamma + M -w) |^4 \leq \int_{B_{2R}} |\grad \eta|^4 + \gamma^{-3}  \int_{B_{2R}} |f|.
\end{equation}
The right hand side of \eqref{jn} has an absolute upper bound, so by the John-Nirenberg Theorem,
there exist absolute constants  $\beta_1, C_* >0$ such that 
\begin{equation} \label{bridge}
\left(\frac {1}{R^4} \int_{B_{2R}}  ( \gamma + M -w)^{\beta_1} \right) \left(\frac {1}{R^4} \int_{B_{2R}} 
 ( \gamma + M -w)^{-\beta_1} \right) \leq C_*.
\end{equation}
Then we can use \eqref{poweriteration} for $\beta= -\beta_1 <0$ and carry out the Moser iteration to 
obtain
\begin{equation} \label{inc}
\left( \frac {1}{R^4} \int_{B_{2R}} (\gamma + M - w)^{-\beta_1} \right)^{- \frac {1}{\beta_1}} \leq \widetilde C^* \inf_{B_R} (\gamma + M - w).
\end{equation}
\eqref{inc} and \eqref{bridge} imply \eqref{whnkMw} for $\beta = \beta_1$. The general case follows from the
case for $\beta = \beta_1$, \eqref{poweriteration} and the Sobolev inequality.
The proof for \eqref{whnkwm} and \eqref{supwm} is done similarly.
\end{proof}
\begin{proof}[Proof of Theorem~\ref{hnkforw}]  \eqref{hnkMw} follows from \eqref{supMw} and \eqref{whnkMw}. 
\eqref{hnkwm}
follows from \eqref{whnkwm} and \eqref{supwm}.
\end{proof}
\begin{proof}[Proof of Theorem~\ref{hnk}]
This will be completed if we can bridge the gap betwen \eqref{upper} and \eqref{lower}. Note that since
$A_{e^{2w}g_0} \in \Gamma_2^+$, $u=e^w$ satisfies $-6 \Delta u = R u^3 >0$, which implies a BMO estimate
for $w = \ln u$ by
\[
0 < \int - \Delta u (u^{-1} \eta^2) =\int |\grad \ln u|^2 \eta^2 + 2 \eta \grad \ln u \cdot \grad \eta,
\]
where $\eta$ is a standard cut-off funtion. So 
\begin{equation} \label{bmo}
\int_{B_R} |\grad w|^2 \leq 4\int_{B_R} |\grad \eta|^2 \lesssim R^2.
\end{equation}
This can also be done as in (2.17) of \cite{CGY2}. Then by the John-Nirenberg Theorem on BMO functions,
there exist absolute constants $\beta_*>0$ and $C_*>0$ such that
\[
\left( \frac{1}{R^4} \int_{B_{2R}} e^{\beta_* w} \right)\left( \frac{1}{R^4} \int_{B_{2R}} e^{-\beta_* w} \right)
\leq C_*.
\]
This estimate, together with \eqref{upper} and \eqref{lower}, conclude the proof of Theorem~\ref{hnk}.
\end{proof}
\begin{proof}[Proof of Theorem~\ref{small}]
The proof consists of two elements:
\begin{enumerate}
\item[i.] There exists absolute constant $\epsilon_0>0$ such that if $w$ is a solution of \eqref{1}
and $\int_{B_{R}} K(x) e^{4w} \leq \epsilon_0$,
then
\begin{equation} \label{gradient}
\int_{B_{R/2}} |\grad w|^4 \leq C,
\end{equation}
where $C$ depends on $ ( \sup_{B_R} K ) \int_{B_R} e^{4w}$.
\item[ii.] Now treat $K(x)e^{4w}$ as $f(x)$. We can verify that, for $p>1$, $R^{4(1-\frac 1p)}||f||_{L^p(B_{R/2})}$
has a bound from above depending only on $ ( \sup_{B_R} K ) \int_{B_R} e^{4w}$,
so that by \eqref{upper} and \eqref{hnkew}, we conclude Theorem~\ref{small}.
\end{enumerate}
To prove \eqref{gradient}, we multiply both sides of \eqref{1} by $\eta^4 (w-\bar w)$, where $\bar w$ is the
average of $w$ over $B_{R}$ and $\eta$ is a standard cut-off function supported in $B_{R}$ such that
for $\rho < R$, $\eta \equiv 1$ on $B_{\rho}$ and $| \grad \eta| \leq 2 (R-\rho)^{-1}$. As before, we have
\begin{equation} \label{last}
\begin{split}
\int_{B_R} \eta^4 |\grad w|^4 &\leq \int_{B_R} \eta^2 |w-\bar w| \left[|\grad w|^2 |\grad \eta|^2 + 
|\grad w|^3 \eta |\grad \eta| \right] + \eta^4 K(x) e^{4w} (w- \bar w) \\
& \leq \left( \int_{B_R} \eta^4 |\grad w|^4 \right)^{1/2} 
\left( \int_{B_R} |w-\bar w|^2 |\grad \eta|^4 \right)^{1/2} + \\
&+ \left( \int_{B_R} \eta^4 |\grad w|^4 \right)^{3/4}\left( \int_{B_R} |w-\bar w|^4 |\grad \eta|^4 \right)^{1/4}
+ \int_{B_R} \eta^4 K(x) e^{4w} (w- \bar w).
\end{split}
\end{equation}
By the BMO estimate \eqref{bmo}, 
\[
\int_{B_R} |w-\bar w|^2 |\grad \eta|^4 \lesssim R^4 (R-\rho)^{-4},
\]
and
\[
\int_{B_R} |w-\bar w|^4 |\grad \eta|^4 \lesssim R^4 (R-\rho)^{-4}.
\]
The last term in \eqref{last} is estimated by Jensen's inequality as
\[
\begin{split}
\int_{B_R} \eta^4 K(x) e^{4w} (w- \bar w) &\leq \left( \int_{B_R}  K(x) e^{4w} \right) \left[
\ln \left( \int_{B_R}  K(x) e^{5w - \bar w} \right) - \ln \left(\int_{B_R}  K(x) e^{4w} \right)\right]\\
&\leq \left( \int_{B_R}  K(x) e^{4w} \right) \left[
\ln \left( \int_{B_R}  K(x) e^{5(w - \bar w)} \right) + 4\bar w - \ln \left(\int_{B_R}  K(x) e^{4w} \right)\right].
\end{split}
\]
By the Moser-Trudinger inequality,
\begin{equation} \label{casep1}
\ln \left( \int_{B_R}  K(x) e^{5(w - \bar w)} \right)\leq c_1 \int_{B_R} |\grad w|^4 +c_2 + \ln (\sup K)+ \ln |B_R|,
\end{equation}
where $c_1, c_2>0$ are absolute constants. 
By Jensen's inequality again, we have
\[
e^{4 \bar w} \leq \frac {1}{|B_R|} \int_{B_R} e^{4w},
\]
which implies 
\[
4 \bar w + \ln |B_R| \leq \ln \left(\int_{B_R} e^{4w}\right).
\]
Substituting all these estimates back into \eqref{last}, we have
\[
\begin{split}
\int_{B_{\rho}} |\grad w|^4 \leq & \left( \int_{B_R}  K(x) e^{4w} \right) \left\{
c_1 \int_{B_R}  |\grad w|^4 +c_2 + \ln \left[\sup K \left(\int_{B_R} e^{4w}\right)/ \left(\int_{B_R}  K(x) e^{4w} \right)\right]\right\}+ \\
&+ C R^4(R-\rho)^{-4} \\
\leq & \epsilon_0 c_1 \int_{B_R}  |\grad w|^4 + \epsilon_0 \left[c_2+ \ln \left(\sup K \int_{B_R} e^{4w}\right)
- \ln \epsilon_0 \right] +  C R^4(R-\rho)^{-4},
\end{split}
\]
if $\int_{B_R}  K(x) e^{4w} \leq \epsilon_0$.  Suppose $\epsilon_0>0$ is chosen so that
 $\epsilon_0 c_1 <1$, then 
from the previously cited Lemma 5.1 in \cite{Gia}, we conclude that, for some absolute constant $\delta >0$, 
\[
\int_{B_{\rho}} |\grad w|^4 \leq \delta \left\{ \epsilon_0 \left[c_2+\ln \left( \sup K \int_{B_R} e^{4w}\right)
-\ln \epsilon_0 \right] + C R^4(R-\rho)^{-4}\right\}.
\]
By taking $\rho = R/2$, we obtain \eqref{gradient}.

We next use \eqref{gradient} to 
estimate $R^{4(1-\frac 1p)} ||f||_{L^p(B_{R/2})}$ with $f = K(x) e^{4w}$ and $p>1$, say,
$p= 5/4$. By the John-Nirenberg
theorem and \eqref{gradient}, we have
\[
\int_{B_{R/2}} e^{5(w-\bar w)} \leq C |B_{R/2}|,
\]
so that
\[
\int_{B_{R/2}} e^{5w} \leq C |B_{R/2}| e^{5\bar w} \leq C |B_{R/2}| \left(\frac {1}{|B_{R/2}|} \int_{B_{R/2}} e^{4w}\right)^{5/4},
\]
and
\begin{equation} \label{casep2}
R^{4(1-\frac 1p)} ||f||_{L^p(B_{R/2})} \leq C R^{4/5}\cdot |B_{R/2}|^{-1/5} \cdot (\sup K) \int_{B_{R/2}} e^{4w}
\leq C (\sup K) \int_{B_{R}} e^{4w}.
\end{equation}
Now we can use \eqref{upper} and \eqref{hnkew} to conclude Theorem 1. When only $||K||_p < \infty$, for
some $p>1$, is assumed, the only changes are in the estimates \eqref{casep1} and \eqref{casep2}. Both can be
handled by first applying the H\"older inequality, then applying the Moser-Trudinger inequality. For instance,
\[
\begin{split}
\ln \left( \int_{B_R}  K(x) e^{5(w -\bar w)} \right) & \leq \ln \left( ||K||_p ||e^{5(w - \bar w)}||_{p'} \right) \\
& \leq \frac{1}{p'}\ln \left( \int_{B_R}   e^{5p'(w - \bar w)} \right) + \ln ||K||_p \\
& \leq c_1 (p')^3 \int_{B_R} |\grad w|^4 +\frac{1}{p'}\left(c_2 + \ln |B_R|\right) + \ln ||K||_p .
\end{split}
\]
Now in choosing $\epsilon_0$, we have to require $\epsilon_0 c_1 (p')^3<1$.
The rest of the modifications are straightforward.
\end{proof} 

\end{document}